\documentclass[12pt,a4paper]{article}
\usepackage[numbers]{natbib}
\usepackage[dvipdfmx]{graphicx}
\usepackage{listings,jvlisting,amsmath,amsfonts}
\usepackage{amsthm}
\usepackage{bm}
\usepackage{color}
\usepackage[subrefformat=parens]{subcaption}
\theoremstyle{definition}
\newtheorem{theorem}{Theorem}
\newtheorem*{theorem*}{Theorem}
\newtheorem{definition}[theorem]{Definition}
\newtheorem*{definition*}{Definition}

\newtheorem*{example*}{Example}

\newcommand{\three}{\mathrm{I}\hspace{-1.2pt}\mathrm{I}\hspace{-1.2pt}\mathrm{I}}
\newcommand{\pedot}[2]{\langle {#1},{#2} \rangle_{\mathrm{pe}}}
\newcommand{\penorm}[1]{\| {#1} \|_{\mathrm{pe}}}
\newcommand{\trans}[1]{\hat{#1}}
\newcommand{\n}{\trans{\mathbf{n}}}
\newcommand{\R}{\mathbb{R}}
\newcommand{\diag}{\mathop{\rm diag}\nolimits}

\begin{document}

\begin{center}
{\Large Laguerre geometry for optimization of gridshell with specified force distribution}

\bigskip\bigskip

Kohei Kabaki$^{\rm a}$, 
Kentaro Hayakawa$^{\rm b}$, 
Makoto Ohsaki$^{\rm a}$

\bigskip
{\it
  ${}^{\rm a}$Department of Architecture and Architectural Engineering, Kyoto
  University, Kyoto-Daigaku Katsura, Nishikyo, Kyoto 615-8540, Japan

  ${}^{\rm b}$Department of Conceptual Design, College of Industrial Technology,
  Nihon University, 1-2-1 Izumi-cho, Narashino, Chiba 275-8575, Japan
  }
\end{center}

\begin{abstract}
For the design of gridshells consisting of continuous beams in two directions and quadrilateral faces to cover a large space of architecture, it is important to arrange each row and column of beams along a planar curve and ensure planar faces for constructability and cross-sectional compatibility at joints.
It is also important that the gridshell is in equilibrium mainly with axial forces against the design loads; i.e., bending deformation should be avoided.
In this study, we first find a continuous shell surface where the principal curvature lines coincide with the principal directions of membrane forces. 
For this purpose, we utilize the L-isothermic surface, which is a kind of membrane O-surface.
Specifically, the generalized Dupin cyclide is used as the reference surface, which has an explicit form of membrane forces with a single arbitrary parameter against normal loads.
Various force distributions are obtained as the parameter is adjusted without changing the surface shape. 
Next, the shell is discretized into a gridshell, and target axial forces are obtained from the section length of the covering region of each node. The axial forces are adjusted by optimizing the cross-sectional radii of the beams with pipe sections to realize the specified target force distribution. 
Since the Laguerre transformation preserves geometric and static properties of the L-isothermic surface, the force can be adjusted by a simple process without re-optimization to obtain an approximate optimal solution of the gridshell after transformation.
The ratio of out-of-plane shear force to the normal load at the node is also evaluated to investigate the effect of deformation on the force distribution.
\end{abstract}

\noindent {\it Keywords}: gridshell, Laguerre geometry,
  generalized Dupin cyclide, \\
  \hspace{2.5cm} cross-sectional optimization

\section{Introduction}
Recently, we have experienced rapid advancement in theory, methodology, and technology of architectural geometry~\cite{Pottmann} for designing continuous or latticed shell structures. 
Among various types of latticed shells, those consisting of continuous beams and quadrilateral faces are called gridshells, and some methods have been proposed for designing the directions of beams of gridshells using geometric properties defined by discrete differential geometry~\cite{Montagne,Pottmann08b,Tellier19}. 
For structural design of a gridshell, it is desirable that it resists the external loads mainly through axial forces without in-plane or out-of-plane bending moments.
Since the mechanical property of a gridshell depends on the surface shape, grid directions, and cross-sections of beams, it is desirable to optimize the shape, geometry, and stiffness simultaneously~\cite{Richardson,Gythiel}.
We can also reduce construction cost and enhance the stiffness by designing the gridshell along planar curves~\cite{Wang}.

One of the simplest approaches to designing a gridshell is to minimize the norm of bending moment or the bending strain energy under the design loads~\cite{Yang}. 
However, in this case, unrealistic distribution may be obtained for the axial forces, and accordingly, the cross-sectional sizes of beam elements.
An alternative approach may be to assign a target desirable distribution of axial forces. However, no good solution is obtained if the target distribution cannot be attained for the specified loading condition or leads to an unrealistic solution.
One possibility for avoiding this inconsistency is to use a known membrane force distribution of a continuous shell and discretize it to derive a target axial force distribution.
However, the critical drawback of this approach is that compatibility in membrane strains is not generally considered in deriving the membrane forces; i.e., only static equilibrium is considered.
Accordingly, the ideal force distribution that can be in equilibrium with the external force without shear cannot be achieved due to the inevitable existence of bending moment and shear force especially near the boundaries.
Nevertheless, it is important to find the target forces equilibrated with the specified loads. An isothermic surface~\cite{Smyth} and a membrane O surface~\cite{Rogers} have favorable properties of resisting uniform pressure loads with membrane forces. 

Recently, M\"obius geometry and Laguerre geometry~\cite{Krasauskas}, which are subsets of Lie sphere geometry~\cite{Bobenko}, have been utilized in surface design of architecture~\cite{Pottmann2007,Mesnil}. 
Cyclide is also a favorable surface, featuring its planar curvature lines, one of which is a circle. 
Schief et al.~\cite{Schief} derived explicit forms of shape and membrane forces of the generalized Dupin cyclide, which has a single arbitrary parameter in the distribution of membrane forces. 
It is also important that the principal directions of membrane forces are aligned to the principal curvature lines~\cite{Pottmann2018,Miki} so that the external loads are transmitted efficiently through axial forces of members after discretization to gridshells.
Kabaki et al.~\cite{Kabaki} utilized the L-minimal generalized Dupin cyclide, which is an L-isothermic L-minimal surface, to find a good target distribution of axial forces of gridshells. 
They optimized the radii of beams to minimize the deviation of axial forces from the target values.
However, Laguerre transformation was not utilized to generate various shapes of gridshells.

In this study, we show that various shapes can be generated through Laguerre transformation preserving the geometric and mechanical properties, and propose a simple adjustment process of section radii to find an approximate optimal solution of the gridshell after transformation.
The ratio of out-of-plane shear force to the normal load is also evaluated to investigate the effect of out-of-plane deformation.

\section{Laguerre geometry and Laguerre transformation}
In this section, basics of Laguerre geometry and Laguerre transformation are summarized for completeness of the paper.

Laguerre geometry is a subgeometry of Lie sphere geometry that manipulates oriented hyperspheres.
In this study, we consider only the geometry in 3-dimensional space. Accordingly, the oriented hypersphere turns out to be the oriented sphere.
Let $\tilde{\mathbf{x}}=(x_1,x_2,x_3)\in\R^3$ and $x_4\in\R$ denote the center and the signed radius of a sphere in the 3-dimensional space. For a plane in the 3-dimensional space, $(x_1,x_2,x_3)$ is the normal vector, and $x_4\in\R$ is its norm. In the following, we mainly consider oriented spheres.

Using the cyclographic model~\cite{Krasauskas}, an oriented sphere is represented by a vector $\mathbf{x} = (x_1,x_2,x_3,x_4)=(\tilde{\mathbf{x}},x_4)\in\R^4$. 
The Euclidean metric in $\R^3$ is extended to the pseudo-Euclidean (\textit{pe}) metric in $\R^4$. The \textit{pe}-inner product of $\mathbf{x}$ and $\mathbf{y}\in\R^4$ is defined as 
\begin{equation}
	\pedot{\mathbf{x}}{\mathbf{y}} = x_1 y_1 + x_2 y_2 + x_3 y_3 - x_4 y_4
\label{pedot}
\end{equation}
Accordingly, the \textit{pe}-norm is defined as
\begin{equation}
	\penorm{\mathbf{x}} = \sqrt{x_1^2 +x_2^2+x_3^2-x_4^2} = \sqrt{\pedot{\mathbf{x}}{\mathbf{x}}}
\end{equation}
The distance between the oriented spheres $\mathbf{x}$ and $\mathbf{y}$ is obtained using the \textit{pe}-norm as
\begin{equation}
        \penorm{\mathbf{x}-\mathbf{y}} = \sqrt{
        \| \tilde{\mathbf{x}}-\tilde{\mathbf{y}} \|^2 - (x_4 - y_4 )^2}
\end{equation}
If $\penorm{\mathbf{x}-\mathbf{y}}=0$, then the two spheres are in \textit{oriented contact}.
Let $\hat{\mathbf{x}}=(x_0,\mathbf{x})^{\top}$ denote a 5-dimensional vector with the weight $x_0$ representing an oriented sphere if $x_0=1$ and a plane if $x_0=0$.
The transformation in Laguerre geometry is called \textit{Laguerre transformation}---a linear transformation in a projective space of the cyclographic model, which is defined as follows~\cite{Krasauskas}:
\begin{definition}{\textbf{Laguerre transformation}\\}
For a vector $\hat{\mathbf{x}}\in\R^5$ with $x_0=1$ or 0, Laguerre transformation is defined as a linear transformation using a pseudo-orthogonal matrix $D\in\R^{4\times4}$, a vector $\mathbf{t}\in\R^4$, and a scalar $\lambda$ as
\begin{equation}
\label{Laguerre}
	\hat{\mathbf{x}} \mapsto \begin{pmatrix}
		1 & 0 \\
		\mathbf{t} & \lambda D \\
	\end{pmatrix} \hat{\mathbf{x}}
\end{equation}
where $D$ satisfies the following relation with $E_\mathrm{pe} = \diag(1,1,1,-1)$:
\begin{equation}
	D^\top E_{\mathrm{pe}} D = E_{\mathrm{pe}}
	\label{constraint_A}
\end{equation}
\end{definition}

\bigskip

It immediately follows from Eq.~(\ref{constraint_A}) that the \textit{pe}-inner product is invariant under Laguerre transformation up to scalar multiplication if $\mathbf{t}=\mathbf{0}$, and the following relation is satisfied:
\begin{equation}
	\pedot{\mathbf{x}}{\mathbf{y}} = \pedot{D \mathbf{x}}{D \mathbf{y}}
\label{keep_dot}
\end{equation}
Let $E\in\R^{4\times4}$ denote the unit matrix.
Equation~(\ref{Laguerre}) corresponds to (i) Scaling if $\lambda\ne 0$, $\mathbf{t}=\mathbf{0}$, $D=E$; (ii) Euclidean transformation and rotation if $\lambda=1$, and $\mathbf{t}$ and $D$ are defined to represent translation and rotation for $(x_1,x_2,x_3)^\top$ with fixed $x_4$; (iii) Offset if $\lambda=1$, $\mathbf{t}=(0,0,0,d)^\top$, $d\ne 0$, and $D=E$; and (iv) \textit{pe}-rotation, for example, the \textit{pe}-rotation about the $x_2 x_3$-plane, which  is defined as follows, with the parameter $\tau_1$ and fixed $x_0$, $x_2$, and $x_3$:
\begin{equation}
	 \begin{pmatrix}
		x_1 \\ x_4 \\
	\end{pmatrix} \mapsto
	 \begin{pmatrix}
		\cosh\tau_1 & \sinh\tau_1 \\
		\sinh\tau_1 & \cosh\tau_1 \\
	\end{pmatrix}
	 \begin{pmatrix}
		x_1 \\ x_4 \\
	\end{pmatrix}
    \label{pe-rotation}
\end{equation}
Similarly, the \textit{pe}-rotations about the $x_3 x_1$- and $x_1 x_2$-planes are defined by the parameters $\tau_2$ and $\tau_3$, respectively.

In the general case,  the vector $\mathbf{t}$ in Eq.~(\ref{Laguerre}) corresponds to translation and offset; the scalar $\lambda$ represents the coefficient of similarity, and $DD$ plays an important role in changing surface shape by transformation. 
Alternatively, the \textit{pe}-norm is an \textit{relative invariant} of weight $\lambda$~\cite{Fels} or an \textit{algebraic invariant} of scalar density of weight 1~\cite{Iri,Fujita} if $\mathbf{t}=\mathbf{0}$.
Note that the distance between the two oriented spheres does not change under Laguerre transformation with $\lambda=0$, and the two contacting spheres remain in contact after transformation with $\mathbf{t}=\mathbf{0}$. 

Since Laguerre geometry in 3-dimensional space manipulates oriented planes and spheres, and cannot directly transform curved surfaces, a point on the surface should be defined as a contact point between the contact elements, i.e., two spheres, or a sphere and a plane.  
See Appendix A for details.
Let $(\xi,\eta)$ and $\mathbf{n}$ denote the curvature-line parameters and the unit normal vector, respectively, of the surface.
The L-isothermic surface, transformed to an L-isothermic surface by Laguerre transformation, is defined as follows:

\begin{definition}{\textbf{L-isothermic surface}}\\The surface is L-isothermic if the third fundamental form $\three$ of the surface is expressed by $(\xi,\eta)$ as
\begin{equation}
	\three = \mathrm{d}\mathbf{n} \cdot \mathrm{d}\mathbf{n} = e^{2\theta} (\mathrm{d}\xi^2 + \mathrm{d}\eta^2)
\label{3_Liso}
\end{equation}
where the dot denotes the inner product in the Euclidean space, and  $\theta$ is a function of $(\xi,\eta)$.
\end{definition}
\bigskip
It is confirmed in Appendix B that an L-isothermic surface is transformed to an L-isothermic surface by Laguerre transformation.
\section{Target membrane forces under normal pressure load}
In this section, we first summarize the results of Schief et al.~\cite{Schief} for completeness of the paper, and show how the target force distribution is defined using a single arbitrary parameter.

Let $T_1$ and $T_2$ denote the normal stresses (normal forces per unit length of section) in the directions of two curvature lines, respectively, parametrized by $\xi$ and $\eta$. In the following, the values corresponding to the curvature lines are denoted by subscripts 1 and 2. The principal curvatures are denoted by $\kappa_1$ and $\kappa_2$, and the lengths of tangent vectors (parametric speeds) along the curvature lines are denoted by $A_1$ and $A_2$. The following equilibrium equations in the tangent and normal directions are satisfied for the uniform normal pressure load $Z$:
\begin{subequations}
\begin{align}
\label{stress_equ_a}
    T_{1\xi} + (\ln A_2)_\xi (T_1-T_2) &= 0 \\
\label{stress_equ_b}
    T_{2\eta}  + (\ln A_1)_\eta  (T_2-T_1) &= 0 \\
\label{stress_equ_c}
\kappa_1 T_1 + \kappa_2 T_2 + Z &= 0
\end{align}
\end{subequations}
where the subscripts $\xi$ and $\eta$ denote the differentiation with respect to $\xi$ and $\eta$, respectively.
By differentiating Eq.~\eqref{stress_equ_c} with respect to $\xi$ and $\eta$, respectively, the following equations are derived:
\begin{equation}
\begin{split}
	T_{1\eta} &= -\left[ \ln (A_1 \kappa_1^2) \right]_\eta T_1 + \frac{\kappa_2}{\kappa_1} \left[ \ln \left(\frac{A_1}{\kappa_2} \right) \right]_\eta T_2 \\
	T_{2\xi} &=\frac{\kappa_1}{\kappa_2} \left[ \ln \left(\frac{A_2}{\kappa_1} \right) \right]_\xi  T_1 - \left[ \ln (A_2 \kappa_2^2) \right]_\xi T_2 \\
\end{split}
\label{membrane-equ2}
\end{equation}
By differentiating Eqs.~\eqref{stress_equ_a} and \eqref{stress_equ_b}, the following condition is derived from the relation $T_{i\xi\eta} = T_{i\eta\xi}$ $ (i = 1,2)$ to be satisfied from invariance with respect to the order of differentiation:
\begin{equation}
	\bar{\mu} T_1 + \bar{\nu} T_2 = 0
\label{geometric_constraint}
\end{equation}
where
\begin{equation}
\begin{split}
	\bar{\mu}	&= \kappa_1 \left\{ \left[ \ln \left(\frac{A_1 \kappa_1}{A_2\kappa_2} \right) \right]_{\xi\eta} + \Upsilon \right\} \\ 
	\bar{\nu}	&= -\kappa_2 \left\{ \left[ \ln \left(\frac{A_1 \kappa_1}{A_2\kappa_2} \right) \right]_{\xi\eta} - \Upsilon \right\} \\
	\Upsilon	&= [\ln (\kappa_1 \kappa_2)]_{\xi\eta} + (\ln A_1)_\eta (\ln \kappa_1)_\xi + (\ln A_2)_\xi (\ln \kappa_2)_\eta - (\ln \kappa_1)_\xi (\ln \kappa_2)_\eta \\
\end{split}
\end{equation}
Equation~\eqref{geometric_constraint} is satisfied independently of the values of $T_1$ and $T_2$ due to the characteristics of the L-isothermic surface satisfying Eq.~\eqref{3_Liso}. Specifically, for canal surfaces where all lines of curvature are planar curves, Schief et al.~\cite{Schief} derived the expressions of $T_1$ and $T_2$ as follows, using an arbitrary parameter $I_0$: 
\begin{equation}
\begin{split}
	T_1 &= -\frac{Z}{2\kappa_2}-\frac{I_0}{\kappa_2 A_2^2} \\
	T_2 &= -\frac{Z}{2\kappa_1} \left[ 1-\frac{(A_1 - A_2)^2}{A_1^2}\right] + \frac{I_0}{\kappa_1 A_1^2} \\
\end{split}
\end{equation}
Note that the first terms in $T_1$ and $T_2$ are at equilibrium with the normal load $Z$, and the second terms correspond to self-equilibrium forces that are proportional to $I_0$.
It is confirmed in Appendices C and D that these solutions are applicable to surfaces after Laguerre transformation and can be analytically calculated.
In the numerical examples, the value of $I_0$ is determined from an additional condition on $T_1$ and $T_2$ to obtain the target distribution of forces.

\section{Cross-sectional optimization of gridshells }

The continuous shell is discretized into a gridshell, and the cross-sectional radii of the beams with pipe sections are optimized~\cite{Kabaki}.
The parameter plane $(\xi,\eta)\in\R^2$ is discretized into a regular grid.
The target axial force and the nodal load are computed from the width and area, respectively, of the covering region of each node. 
Structural response under normal loads is computed using standard finite element analysis with beam elements, where the beams have pipe sections.

Let $N_i$ and $N_i^*$ $(i=1,\dots,n)$ denote the axial force and its target value of beam element $i$, respectively, where $n$ is the number of elements.
The ratio of thickness to external radius of the pipe section is fixed, and the radius is taken as the design variable.
The continuously located beams in the same plane are assumed to be in the same group with the same radius.
The radius of the beams in group $j$ is denoted by $R_j$, and the radii of all groups are combined to the design variable vector $\mathbf{R}$ with its upper and lower bounds $\mathbf{R}^{\mathrm{U}}$ and $\mathbf{R}^{\mathrm{L}}$, respectively. 

The optimization problem is formulated as
\begin{subequations}
\begin{alignat}{2}
 & {\rm Minimize} & \ \  & F(\mathbf{R}) = \sum_{i=1}^{n} (N_i(\mathbf{R})-N_i^*)^2 \\
 & {\rm subject\ to} & \ \  & \mathbf{R}^{\mathrm{L}} \le \mathbf{R} \le \mathbf{R}^{\mathrm{U}}
\end{alignat}
\end{subequations}
which is solved using a nonlinear programming library in the numerical examples.
Note that the objective function of this optimization problem cannot be reduced to 0, if the ideal force distribution that can be in equilibrium with the external force without shear force is assigned as the target force, due to the inevitable existence of bending moment and out-of-plane shear force especially near the boundaries.

Since the properties of Laguerre geometry are preserved after Laguerre transformation, the optimal solution before transformation can be expected to be a good initial solution for optimizing the cross-sectional radii after transformation.
The optimal cross-sectional area of member $i$ before transformation is denoted by $\hat{S}_i$.
Let $\hat{\mathbf{R}}$ denote the optimal solution before transformation.
The axial force of member $i$ of the gridshell obtained by finite element analysis for the solution $\hat{\mathbf{R}}$ with the shape after transformation is denoted by $N_i^1(\hat{\mathbf{R}})$.
The cross-sectional area is adjusted by the following simple stress ratio formula:
\begin{equation}
    \label{adjust}
    S_i =  \hat{S}_i \frac{N_i^1(\hat{\mathbf{R}})}{N_i^*}
\end{equation}
where $N_i^*$ is the target value evaluated for the shape after transformation. The radius $R_i$ is updated from the area $S_i$. 

\section{Numerical examples}
The shell surface is obtained by cutting a region from the generalized Dupin cyclide, which is a canal surface with a sphere of varying radius.
The gridshell is then generated by discretizing the surface into a grid.

\begin{figure}
\centering
\begin{minipage}[b]{0.49\columnwidth}
    \centering
    \includegraphics[width=0.9\columnwidth]{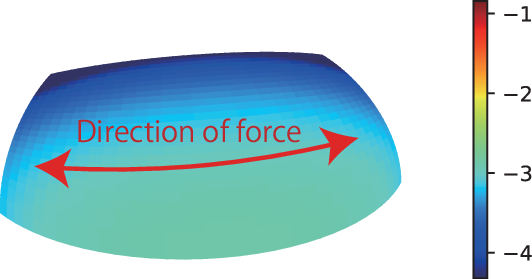}
    \subcaption{}
    \label{NL_target_T1}
\end{minipage}
\begin{minipage}[b]{0.49\columnwidth}
    \centering
    \includegraphics[width=0.9\columnwidth]{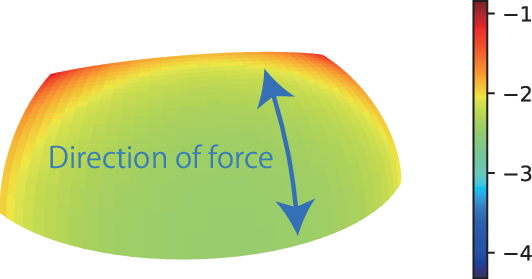}
    \subcaption{}
    \label{NL_target_T2}
\end{minipage}
\caption{Target membrane force distribution (N/mm) before transformation; (a) Target $T_1$ in the direction of generating line in red arrow, (b) Target $T_2$ in the direction of circle in blue arrow.}
\label{target_membrane}
\end{figure}

The surface is expressed by Schief et al.~\cite{Schief} with respect to the parameters $\alpha$ and $u$, which are replaced by $\xi$ and $\eta$, respectively, as 
\begin{equation}
    \mathbf{r} = \frac{\cos\xi}{2(a_0-c_0 \cos\xi\cos\eta)}
    \begin{pmatrix}
        \cos\xi \\ a_0 \sin\xi \\ \cos\xi\sin\eta \\
    \end{pmatrix}
    + \frac{1}{4}
    \begin{pmatrix}
        -a_0 \\ 2\xi \\ 0 \\
    \end{pmatrix}
    \label{Dupin}
\end{equation}
where $a_0$ and $c_0$ are the parameters satisfying $a_0^2-c_0^2=1$.  
Consider a region defined by $-0.1\pi\le\xi\le 0.1\pi$ and $0\le\eta\le 0.15$, for a surface corresponding to $(a_0,c_0)=(2.144,1.896)$.
As an example we generate a gridshell so that $T_1 A_2=T_2 A_1$ is satisfied at the center $(\xi,\eta)=(0,0.075)$, where $\kappa_1$, $\kappa_2$ and $A_2$ are computed numerically by finite difference approximation, and $A_1$ is obtained analytically differentiating $\mathbf{r}$ in Eq.~\eqref{Dupin} with respect to $\xi$.
The value of $I_0$ in this case is $-76850$ N for $Z= -0.0005 \; \mathrm{N/mm}^2$. 
The target distributions of $T_1$ and $T_2$ are plotted in Fig.~\ref{target_membrane}.
Note that $T_1$ corresponds to $\eta$ in the direction of generating curve of the canal surface (the red arrow in Fig.~\ref{target_membrane}(a)), while $T_2$ corresponds to $\xi$ in the direction of circle (blue arrow in Fig.~\ref{target_membrane}(b)), 

\begin{figure}
\centering
    \includegraphics[width=0.6\columnwidth]{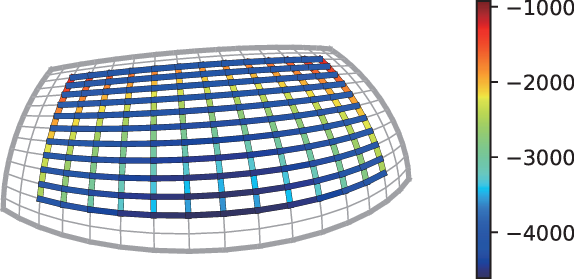}
    \caption{Target axial force (N) of gridshell excluding two rows of members along boundary before transformation.}
    \label{NL_target_N}
\end{figure}

The surface is discretized into a gridshell so that each beam covers the section of the region corresponding to the grid size determined by uniformly dividing the parameters $\xi\in[-0.1\pi,0.1\pi]$ and $\eta\in[0,0.15]$ into 14 and 16 grids, respectively.
The obtained target force distribution is plotted in Fig.~\ref{NL_target_N}, where the members in $\eta$ direction have large compressive forces.
Note that the members in the two rows along the boundary are not considered in computing the norm of force difference in the objective function, because it is difficult to adjust the axial forces of members near the boundary by modifying only the axial stiffness.
Accordingly, the radii of the members in two rows along the boundary are not included in the design variables.

\begin{figure}
\centering
\begin{minipage}[b]{0.49\columnwidth}
    \centering
    \includegraphics[width=0.9\columnwidth]{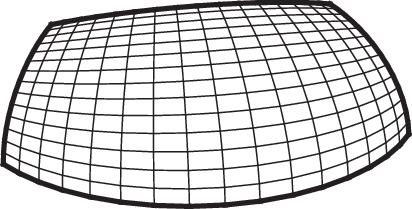}
    \subcaption{}
\end{minipage}
\begin{minipage}[b]{0.49\columnwidth}
    \centering
    \includegraphics[width=0.9\columnwidth]{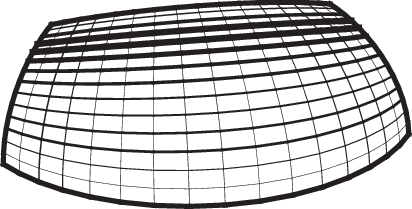}
    \subcaption{}
\end{minipage}
\caption{Cross-sectional radius of solutions before transformation. The line width is proportional to beam radius; (a) Initial solution, (b) Optimal solution. }
    \label{cross-sectiona_radius}
\end{figure}

\begin{figure}
\centering
\begin{minipage}[b]{0.49\columnwidth}
    \centering
    \includegraphics[width=0.9\columnwidth]{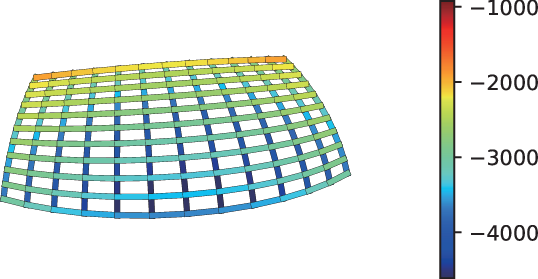}
    \subcaption{Initial}
\end{minipage}
\begin{minipage}[b]{0.49\columnwidth}
    \centering
    \includegraphics[width=0.9\columnwidth]{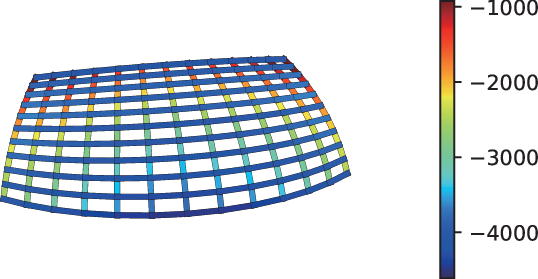}
    \subcaption{Optimal}
\end{minipage}
\caption{Axial force (N) of solutions before transformation; (a) Initial solution, (b) Optimal solution.}
    \label{axial_force}
\end{figure}

The boundary nodes are pin supported, and the linear elastic analysis is carried out using Python 3.11.3, where the ratio of thisckness to radius is 0.1, all members have the same Young's modulus, and shear deformaion is not considered.
The gridshell is scaled by a factor of 10000 mm to investigate the structural property in a realistic size.
The optimization problem is solved to optimize the cross-sectional radii of gridshells.
The nonlinear programming program Nelder-Mead in the SciPy library Ver. 1.24.3 is used for optimization.
The initial and optimal solutions are shown in Fig.~\ref{cross-sectiona_radius}(a) and (b), respectively, where the width of each member is proportional to its radius.
In this case, the members in the upper region became stiffer after optimization.
The axial force distributions of initial and optimal solutions are shown in Fig.~\ref{axial_force}(a) and (b), respectively.
We can see from Figs.~\ref{NL_target_N} and \ref{axial_force}(b) that the force distribution close to the target has been obtained by optimization.

Table~\ref{opt_original} shows the maximum and mean differences of axial forces from the target values, which are reduced to 33\% and 18\%, respectively, after optimization.
Let $Q_i$ denote the sum of out-of-plane shear force from the beams connected to node $i$.
The ratio of $|Q_i|$ to the normal load at node $i$ is computed at all nodes except boundary nodes and their mean values are listed in Table~\ref{opt_original}, which shows that the effect of out-of-plane shear force slightly increases as a result of optimization.
However, the optimal solution mainly resists the normal loads through axial forces.

\begin{table}[tbp]
    \centering
  \caption{Optimization result of gridshell before transformation: Deviation from the target axial forces and shear force ratio.}
  \label{opt_original}
    \begin{tabular}{lccc}
        \hline
      & Max. $|N_i-N_i^*|$ & Mean $|N_i-N_i^*|$ & Mean shear/load \\
    \hline \hline
    Initial     & 2135      & 1323  & 0.124\\
    Optimal     & 707       & 239   & 0.146\\
    \hline
  \end{tabular}
\end{table}

\begin{figure}
\centering
\begin{minipage}[b]{0.32\columnwidth}
    \centering
    \includegraphics[width=0.95\columnwidth]{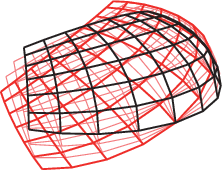}
    \subcaption{}
\end{minipage}
\begin{minipage}[b]{0.32\columnwidth}
    \centering
    \includegraphics[width=0.95\columnwidth]{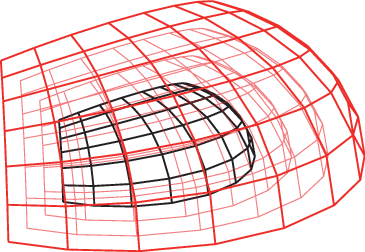}
    \subcaption{}
\end{minipage}
\begin{minipage}[b]{0.32\columnwidth}
    \centering
    \includegraphics[width=0.95\columnwidth]{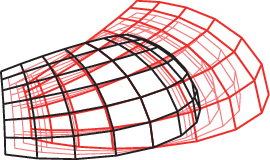}
    \subcaption{}
\end{minipage}
\caption{Examples of Laguerre transformation; (a) Euclidean rotation, (b) Offset, (c) \textit{pe}-rotation in projective space.}
    \label{ex_Laguerre}
\end{figure}

\begin{figure}
\centering
\begin{minipage}[b]{0.49\columnwidth}
    \centering
    \includegraphics[width=0.9\columnwidth]{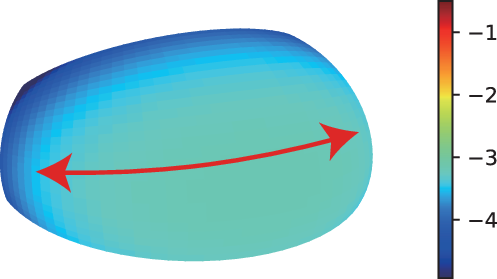}
    \subcaption{Target $T_1$}
    \label{tr_target_T1}
\end{minipage}
\begin{minipage}[b]{0.49\columnwidth}
    \centering
    \includegraphics[width=0.9\columnwidth]{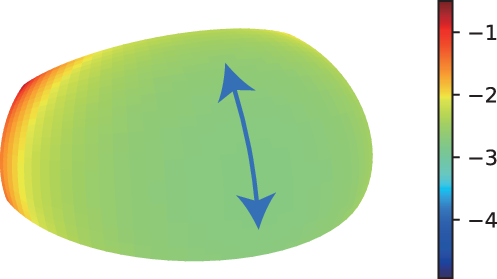}
    \subcaption{Target $T_2$}
    \label{tr_target_T2}
\end{minipage}
\caption{Target membrane force distribution (N/mm) after transformation.}
\label{membrane_force_trans}
\end{figure}

\begin{figure}
\centering
    \includegraphics[width=0.6\columnwidth]{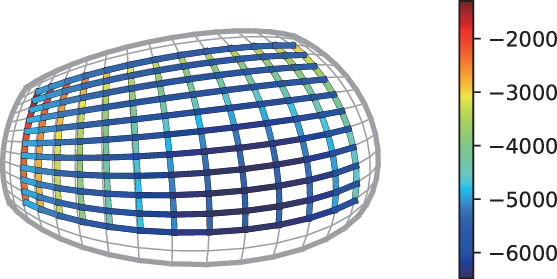}
    \caption{Target axial force (N) of gridshell excluding two rows of members along boundary after transformation.}
    \label{target_trans}
\end{figure}

Among various types of Laguerre transformation in Fig.~\ref{ex_Laguerre}, rotation is a simple transformation that exists even in Euclidean geometry. Offset is a unique type of transformation that blows up the surface in a constant distance in the direction of surface normal.
Here we apply the \textit{pe}-rotation in the projective space using hyperbolic functions, which is also a unique type of transformation of Laguerre geometry as formulated in Eq.~\eqref{pe-rotation}, with the hyperbolic angle parameters $\tau_1 = 0.1$, $\tau_1 = 0.2$, and $\tau_3 = 0.3$. This transformation leads to a surface without symmetry. When the four corner nodes are located on a plane, their coordinates (mm) become: $(0, 0)$, $(2.343\times 10^4, 0)$, $(9.418\times 10^2, 9.234\times 10^3)$, and $(1.866\times 10^4, 1.352\times 10^4)$, and the rise of the surface is 10629 mm.
The target distribution of membrane forces is obtained, as plotted in Fig.~\ref{membrane_force_trans}, so that $T_1 A_2=T_2 A_1$ is satisfied at the center $(\xi,\eta)=(0,0.075)$. 
Target axial force after discretization is shown in Fig.~\ref{target_trans}.

\begin{figure}
\centering
\begin{minipage}[b]{0.49\columnwidth}
    \centering
    \includegraphics[width=0.8\columnwidth]{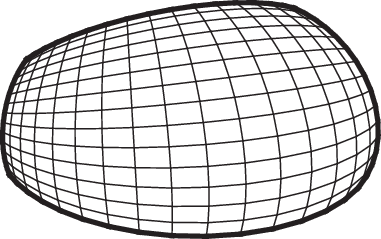}
    \subcaption{}
\end{minipage}
\begin{minipage}[b]{0.49\columnwidth}
    \centering
    \includegraphics[width=0.8\columnwidth]{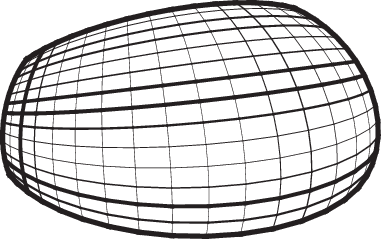}
    \subcaption{}
\end{minipage}
\\
\centering
\begin{minipage}[b]{0.49\columnwidth}
    \centering
    \includegraphics[width=0.8\columnwidth]{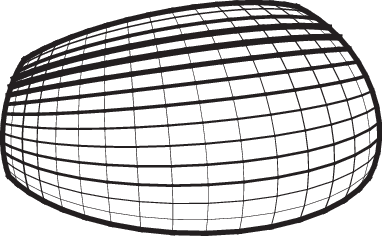}
    \subcaption{}
\end{minipage}
\begin{minipage}[b]{0.49\columnwidth}
    \centering
    \includegraphics[width=0.8\columnwidth]{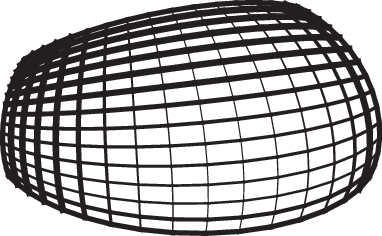}
    \subcaption{}
    \label{tr_after}
\end{minipage}
 \caption{Cross-sectional radius after transformation;  (a) Initial, (b) Optimal solution by Nelder-Mead algorithm, (c) Before adjustment (radii equal to the optimal solution before transformation), (d) After adjustment using Eq.~\eqref{adjust}.}
   \label{cross-section_trans}
\end{figure}

\begin{figure}
\centering
\begin{minipage}[b]{0.49\columnwidth}
    \centering
    \includegraphics[width=0.9\columnwidth]{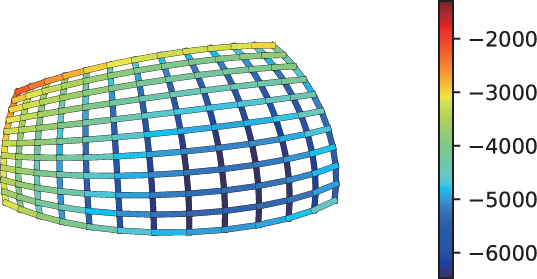}
    \subcaption{}
\end{minipage}
\begin{minipage}[b]{0.49\columnwidth}
    \centering
    \includegraphics[width=0.9\columnwidth]{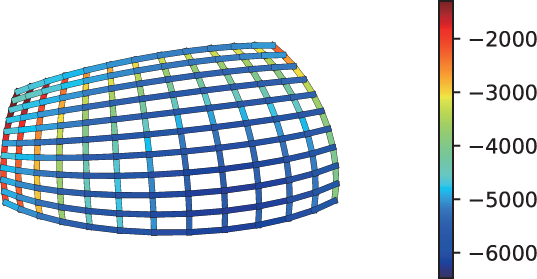}
    \subcaption{}
\end{minipage}
\\
\begin{minipage}[b]{0.49\columnwidth}
    \centering
    \includegraphics[width=0.9\columnwidth]{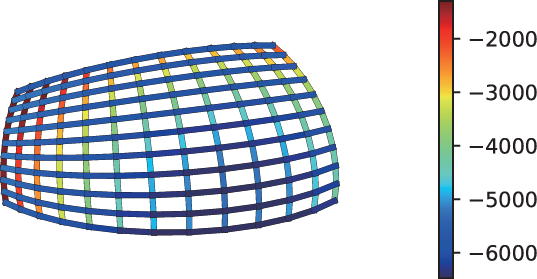}
    \subcaption{}
\end{minipage}
\begin{minipage}[b]{0.49\columnwidth}
    \centering
    \includegraphics[width=0.9\columnwidth]{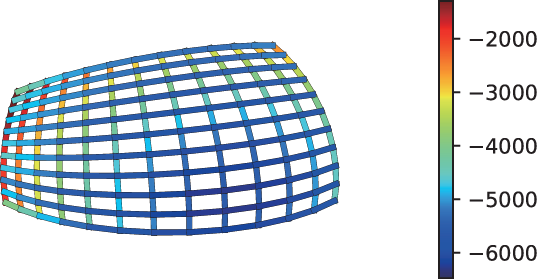}
    \subcaption{}
\end{minipage}
\caption{Axial forces after transformation; (a) Initial, (b) Optimal solution by Nelder-Mead algorithm, (c) Before adjustment (radii equal to the optimal solution before transformation), (d) After adjustment using Eq.~\eqref{adjust}.}
    \label{axial_force_trans}
\end{figure}

Optimization is carried out also for the transformed surface.
The initial and optimal cross-sectional radii are plotted in Fig.~\ref{cross-section_trans}(a) and (b), respectively.
On the other hand, the optimal solution before transformation is directly applied as Figs.~\ref{cross-section_trans}(c), and adjustment in Eq.~\eqref{adjust} is carried out to obtain the solution in Fig.~\ref{cross-section_trans}(d). Axial forces of these solutions are shown in Fig.~\ref{axial_force_trans}.
It is seen from Table~\ref{optimization_result_trans} that the maximum values of errors of the axial forces from the target values is reduced from the initial solution to 114\% of the optimal value by the simple adjustment process, and the ratio to the optimal value is 191\% even before adjustment, where the optimal solution before adjustment is simply used.
The mean value after adjustment is very close to the optimal value.
The ratio of out-of-plane shear force to the normal load is small enough also for the transformed surface.

\begin{table}[t]
    \centering
  \caption{Optimization results after transformation: Deviation from the target axial forces and shear force ratio.}
  \label{optimization_result_trans}
    \begin{tabular}{lccc}
        \hline
      & Max. $|N_i-N_i^*|$ & Mean $|N_i-N_i^*|$ & Mean shear/load \\
    \hline \hline
    Initial     & 2266      & 1407  & 0.073\\
    Optimal     & 821       & 262   & 0.074\\
    Before adjustment     & 1565       & 474   & 0.076\\
    After adjustment     & 941       & 266   & 0.078\\
    \hline
  \end{tabular}
\end{table}

\newpage
\section{Conclusions}

A method has been proposed to minimize the error of axial forces from their target values of gridshells subjected to uniform normal loads.
The generalized Dupin cyclide, which is categorized as an L-isothermic surface, is used as the reference surface, and a surface shape without symmetry is generated by \textit{pe}-rotation in projective space of Laguerre geometry, which is a unique type of Laguerre transformation that has not been utilized in the design of gridshells in architecture.
The benefit of using Laguerre transformation is that the L-isothermic property is preserved by transformation, and explicit forms of membrane forces against normal loads are available also for the transformed surface.
It has been shown in the numerical examples that the member forces can be adjusted, without re-optimization, by a simple updating process of cross-sectional area that is similar to a stress ratio algorithm.
Effect of out-of-plane shear force has been also investigated to show that the optimal solution mainly resists the normal loads through axial forces even considering deformation of the gridshell.

\section*{Acknowledgment}
This research is partly supported by JST Crest JP-MJCR1911.
The authors would like to thank Prof. Yohei Yokosuka at Kagoshima University and Prof. Yoshiki Jikumaru at Toyo University for their valuable comments at the preliminary research stage.

\bibliographystyle{elsarticle-num} 
\bibliography{bib.bib}

\section*{Appendix}
In the following, the properties of Laguerre transformation are presented for the case with $\mathbf{t}=\mathbf{0}$ and $\lambda=1$ because these are obviously satisfied for translation and scaling.

\subsection*{A. Laguerre transformation of curved surface}

The family of oriented spheres $\mathbf{s}(\alpha)$ that contact with each other at point $\mathbf{p}$ on the surface is expressed using the unit normal vector $\mathbf{n}$ of the surface as
\begin{equation}
	\mathbf{s}(\alpha) = 
    \begin{pmatrix}
		\mathbf{p} \\
		0 \\
	\end{pmatrix} + \alpha \begin{pmatrix}
		\mathbf{n} \\
		-1 \\
	\end{pmatrix}
\label{init_surf}
\end{equation}
where $\alpha\in\R$ is an arbitrary real value representing the radius of the sphere.
From the property of Laguerre transformation, which transforms an oriented sphere including oriented plane and a point to an oriented sphere, these spheres contact also after Laguerre transformation. Thus, the surface is defined as the set of points after transformation from the points on the original surface. Laguerre transformation of $\mathbf{s}(\alpha)$ in Eq.~\eqref{init_surf} leads to
\begin{equation}
	D \mathbf{s}(\alpha)
    = D \begin{pmatrix}
		\mathbf{p} \\
		0 \\
	\end{pmatrix} + \alpha D \begin{pmatrix}
		\mathbf{n} \\
		-1 \\
	\end{pmatrix}
\label{lagtra}
\end{equation}
Note that the center of sphere moves along a line called pencil, which is in the direction of the surface normal also after transformation. 
Denote the vector $D(\mathbf{n}^\top,-1)^\top$ in Eq.~(\ref{lagtra}) as $(\n^\top,-a)^\top$. From Eq.~(\ref{constraint_A}) and $\|\mathbf{n}\|=1$, the following relation is satisfied:
\begin{equation}
	\penorm{(\n^\top,-a)} = \penorm{(\mathbf{n}^\top,-1)} = 0
\end{equation}
which leads to $a=\pm\|\n\|$, and contact condition is preserved. Hence, the following relation hols:
\begin{equation}
	D \mathbf{s}(\alpha)
    = D \begin{pmatrix}
		\mathbf{p} \\
		0 \\
	\end{pmatrix} + \alpha \begin{pmatrix}
		\n \\
		\pm\|\n\| \\
	\end{pmatrix}
\label{lagtra2}
\end{equation}

\subsection*{B. Transformation of L-isothermic surface}

In the following, the subscripts $\xi$ and $\eta$ indicate partial differentiation with respect to $\xi$ and $\eta$, respectively. 
For the L-isothermic surface, the following relations hold:
\begin{equation}
	\mathbf{n}_\xi \cdot \mathbf{n}_\xi = \mathbf{n}_\eta \cdot \mathbf{n}_\eta = e^{2\theta},\quad
	\mathbf{n}_\xi \cdot \mathbf{n}_\eta = 0
\end{equation}
The third fundamental form is expressed, as follows, using the normal vector $\n$ defined in Eq.~\eqref{lagtra2}:
\begin{equation}
\begin{split}
	\trans{\three} &= \frac{d\n}{\|\n\|} \cdot \frac{d\n}{\|\n\|} \\
	&= \left(\frac{\n}{\|\n\|} \right)_\xi \cdot \left(\frac{\n}{\|\n\|} \right)_\xi d\xi^2 +2 \left(\frac{\n}{\|\n\|} \right)_\xi \cdot \left(\frac{\n}{\|\n\|} \right)_\eta d\xi d\eta + \left(\frac{\n}{\|\n\|} \right)_\eta \cdot \left(\frac{\n}{\|\n\|} \right)_\eta d\eta^2
\end{split}
\label{trans_three}
\end{equation}
Each term in Eq.~\ref{trans_three} is rewritten as
\begin{equation}
\begin{split}
	\left(\frac{\n}{\|\n\|} \right)_\xi \cdot \left(\frac{\n}{\|\n\|} \right)_\xi &= \frac{1}{\|\n\|^4} (\n_\xi \|\n\| - \n \|\n\|_\xi)\cdot (\n_\xi \|\n\| - \n \|\n\|_\xi) \\
	&= \frac{1}{\|\n\|^4} ( \|\n\|^2 \n_\xi \cdot \n_\xi - 2 \|\n\| \|\n\|_\xi \n \cdot \n_\xi + \|\n\|_\xi^2 \n \cdot \n)  \\ 
	\left(\frac{\n}{\|\n\|} \right)_\xi \cdot \left(\frac{\n}{\|\n\|} \right)_\eta &= \frac{1}{\|\n\|^4} (\n_\xi \|\n\| - \n \|\n\|_\xi)\cdot (\n_\eta \|\n\| - \n \|\n\|_\eta) \\
	&= \frac{1}{\|\n\|^4} (\|\n\|^2 \n_\xi \cdot \n_\eta - \|\n\| \|\n\|_\xi \n \cdot \n_\eta - \|\n\| \|\n\|_\eta \n \cdot \n_\xi + \|\n\|_\xi \|\n\|_\eta \n \cdot \n )  \\ 
	\left(\frac{\n}{\|\n\|} \right)_\eta \cdot \left(\frac{\n}{\|\n\|} \right)_\eta &= \frac{1}{\|\n\|^4} (\n_\eta \|\n\| - \n \|\n\|_\eta)\cdot (\n_\eta \|\n\| - \n \|\n\|_\eta) \\
	&= \frac{1}{\|\n\|^4} ( \|\n\|^2 \n_\eta \cdot \n_\eta - 2 \|\n\| \|\n\|_\eta \n \cdot \n_\eta + \|\n\|_\eta^2 \n \cdot \n)  \\ 
\end{split}
\label{three_coef}
\end{equation}
The following equation is satisfied from Eqs.~\eqref{pedot} and \eqref{keep_dot}:
\begin{equation}
\begin{split}
	\n_\xi \cdot \n_\xi 	&= \pedot{\begin{pmatrix} \n \\ -\|\n\| \\ \end{pmatrix}_\xi}{\begin{pmatrix} \n \\ -\|\n\| \\ \end{pmatrix}_\xi} + \|\n\|_\xi^2 \\
				&= \pedot{D\begin{pmatrix} \mathbf{n} \\ -1 \\ \end{pmatrix}_\xi}{D \begin{pmatrix} \mathbf{n} \\ -1 \\ \end{pmatrix}_\xi} + \|\n\|_\xi^2 \\
				&= \pedot{\begin{pmatrix} \mathbf{n} \\ -1 \\ \end{pmatrix}_\xi}{\begin{pmatrix} \mathbf{n} \\ -1 \\ \end{pmatrix}_\xi} + \|\n\|_\xi^2 \\
				&= \mathbf{n}_\xi \cdot \mathbf{n}_\xi + \|\n\|_\xi^2 \\
				&= e^{2\theta} + \|\n\|_\xi^2 \\
\end{split}
\end{equation}
\begin{equation}
\begin{split}
	\n \cdot \n_\xi 	&= \pedot{\begin{pmatrix} \n \\ -\|\n\| \\ \end{pmatrix}}{\begin{pmatrix} \n \\ -\|\n\| \\ \end{pmatrix}_\xi} + \|\n\| \|\n\|_\xi \\
				&= \pedot{D\begin{pmatrix} \mathbf{n} \\ -1 \\ \end{pmatrix}}{D \begin{pmatrix} \mathbf{n} \\ -1 \\ \end{pmatrix}_\xi} + \|\n\| \|\n\|_\xi \\
				&= \pedot{\begin{pmatrix} \mathbf{n} \\ -1 \\ \end{pmatrix}}{\begin{pmatrix} \mathbf{n} \\ -1 \\ \end{pmatrix}_\xi} + \|\n\| \|\n\|_\xi \\
				&= \mathbf{n} \cdot \mathbf{n}_\xi + \|\n\| \|\n\|_\xi \\
\end{split}
\end{equation}
\begin{equation}
\begin{split}
	\n_\xi \cdot \n_\eta 	&= \pedot{\begin{pmatrix} \n \\ -\|\n\| \\ \end{pmatrix}_\xi}{\begin{pmatrix} \n \\ -\|\n\| \\ \end{pmatrix}_\eta} + \|\n\|_\xi \|\n\|_\eta \\
				&= \pedot{D\begin{pmatrix} \mathbf{n} \\ -1 \\ \end{pmatrix}_\xi}{D \begin{pmatrix} \mathbf{n} \\ -1 \\ \end{pmatrix}_\eta} + \|\n\|_\xi \|\n\|_\eta \\
				&= \pedot{\begin{pmatrix} \mathbf{n} \\ -1 \\ \end{pmatrix}_\xi}{\begin{pmatrix} \mathbf{n} \\ -1 \\ \end{pmatrix}_\eta} + \|\n\|_\xi \|\n\|_\eta \\
				&= \mathbf{n}_\xi \cdot \mathbf{n}_\eta + \|\n\|_\xi \|\n\|_\eta \\
				&= \|\n\|_\xi \|\n\|_\eta \\
\end{split}
\end{equation}
\begin{equation}
	\n \cdot \n 		= \|\n\| ^2
\end{equation}
Since $\mathbf{n}$ is the unit normal vector, $\mathbf{n} \cdot \mathbf{n}_\xi = \mathbf{n} \cdot \mathbf{n}_\eta = 0$ is satisfied, which is incorporated into Eq.~\eqref{three_coef} to obtain
\begin{equation}
\begin{split}
	\left(\frac{\n}{\|\n\|} \right)_\xi \cdot \left(\frac{\n}{\|\n\|} \right)_\xi &= \frac{1}{\|\n\|^4} ( \|\n\|^2 (e^{2\theta} + \|\n\|_\xi^2) - 2 \|\n\| \|\n\|_\xi \|\n\| \|\n\|_\xi  + \|\n\|_\xi^2 \|\n\|^2)  \\
	&= \frac{1}{\|\n\|^2} e^{2\theta} \\ 
	\left(\frac{\n}{\|\n\|} \right)_\xi \cdot \left(\frac{\n}{\|\n\|} \right)_\eta &= \frac{1}{\|\n\|^4} (\|\n\|^2 \|\n\|_\xi \|\n\|_\eta - \|\n\| \|\n\|_\xi \|\n\| \|\n\|_\eta  \\
    & \hspace{2cm}- \|\n\| \|\n\|_\eta \|\n\| \|\n\|_\xi + \|\n\|_\xi \|\n\|_\eta \|\n\|^2 )  \\ 
	&= 0 \\
	\left(\frac{\n}{\|\n\|} \right)_\eta \cdot \left(\frac{\n}{\|\n\|} \right)_\eta &= \frac{1}{\|\n\|^4} ( \|\n\|^2 (e^{2\theta} + \|\n\|_\eta^2) - 2 \|\n\| \|\n\|_\eta \|\n\| \|\n\|_\eta  + \|\n\|_\eta^2 \|\n\|^2)  \\
	&= \frac{1}{\|\n\|^2} e^{2\theta}
\end{split}
\end{equation}
By incorporating a new variable $\trans{\theta}=\theta - \log \|\n\|$ into Eq.~\eqref{trans_three}, the third fundamental form of the transformed surface becomes
\begin{equation}
	\trans{\three} = e^{2\trans{\theta}}(d\xi^2 + d\eta^2)
\end{equation}
which is in the same form as the original surface.

\subsection*{C. Static property of L-isothermic surface}
We consider a generalized Dupin cyclide \cite{Schief} that is a kind of L-minimal surface which is also an L-isothermic surface.
Dupin cyclide is a canal surface generate as an envelope of a sphere with variable radius translating along a curve.
The two variables $T_1$ and $T_2$ may be generally obtained from the third equation of Eq.~\eqref{stress_equ_c} and Eq.~\eqref{geometric_constraint}. However, for the L-isothermic surface, the following relation is satisfied
\begin{equation}
	\kappa_1 A_1 = \kappa_2 A_2 = e^\theta
\end{equation}
which leads to
\begin{equation}
  \label{iso2}
	\ln \left(\frac{A_1 \kappa_1}{A_2 \kappa_2} \right) = 0
\end{equation}
By incorporating Eq.~\eqref{iso2} into Eq.~\eqref{geometric_constraint}, we obtain $(\kappa_1 T_1 +\kappa_2 T_2) \Upsilon = 0$. Furthermore, $\Upsilon = 0$ is derived from Eq.~\eqref{stress_equ_c} and $Z \ne 0$. Therefore, the invariance condition \eqref{geometric_constraint} with respect to the order of differentiation is always satisfied for L-isothermic surface; accordingly, two force components $T_1$ and $T_2$ should satisfy a single equation, and an arbitrary parameter remains.
In the following examples, we consider a canal surface where all curvature lines are planar curves, and it is a class of L-isothermic surface. For a canal surface, $T_1$ and $T_2$ satisfying Eqs.~\eqref{stress_equ_a} and \eqref{stress_equ_b} are explicitly obtained, as follows, with an arbitrary parameter $I_0$:
\begin{equation}
\begin{split}
	T_1 &= -\frac{Z}{2\kappa_2}-\frac{I_0}{\kappa_2 A_2^2} \\
	T_2 &= -\frac{Z}{2\kappa_1} \left[ 1-\frac{(A_1 - A_2)^2}{A_1^2}\right] + \frac{I_0}{\kappa_1 A_1^2} \\
\end{split}
\label{origin_tention_equ}
\end{equation}
As discussed in Appendix B, an L-isothermic surface moves to an L-isothermic surface by Laguerre transformation. For a canal surface, this is obvious because it is generated by translating a sphere of variable radius along a planar curve. Therefore, the normal stresses $\trans{T_1}$ and $\trans{T_2}$ after transformation is written, as follows, using the geometrical properties $\trans{\kappa}_1$, $\trans{A}_1$, $\trans{\kappa}_2$, and $\trans{A}_2$, in the same manner as \eqref{origin_tention_equ}:
\begin{equation}
\begin{split}
	\trans{T}_1 &= -\frac{Z}{2 \trans{\kappa}_2} - \frac{\trans{I}_0}{\trans{\kappa}_2 \trans{A}_2^2} \\ 
	\trans{T}_2 &= -\frac{Z}{2 \trans{\kappa}_1} \left[ 1-\frac{(\trans{A}_1-\trans{A}_2)^2}{\trans{A}_1^2} \right] + \frac{\trans{I}_0}{\trans{\kappa}_1 \trans{A}_1^2} \\
\end{split}
\label{trans_tention_equ}
\end{equation}
Therefore, it is necessary to derive expressions of  $\trans{\kappa}_1$, $\trans{A}_1$, $\trans{\kappa}_2$, and $\trans{A}_2$ after transformation.
\subsection*{D. Geometrical properties after transformation}
Consider two oriented spheres contacting at point $\mathbf{p}$ on the surface, where the radii of spheres are equal to $1/\kappa_1$ and $1/\kappa_2$, respectively. The plane contacting the surface at $\mathbf{p}$ is also considered. These spheres and plane are defined by 5-dimensional vector as
\begin{equation}
	\text{Plane:}
	\begin{pmatrix}
		0 \\ \mathbf{n} \\ -1 \\
	\end{pmatrix} \quad
	\text{Sphere 1:}
	\begin{pmatrix}
		1 \\ \mathbf{p} - \frac{1}{\kappa_1} \mathbf{n} \\ \frac{1}{\kappa_1} \\
	\end{pmatrix} \quad
	\text{Sphere 2:}
	\begin{pmatrix}
		1 \\ \mathbf{p} - \frac{1}{\kappa_2} \mathbf{n} \\ \frac{1}{\kappa_2} \\
	\end{pmatrix}
\end{equation}
Application of Laguerre transformation to the plane and spheres is written as
\begin{equation}
	\begin{pmatrix}
		0 \\ \mathbf{n} \\ -1 \\
	\end{pmatrix} \mapsto
	\begin{pmatrix}
		0 \\ \n \\ -\|\n\| \\
	\end{pmatrix} ,\quad	
	\begin{pmatrix}
		1 \\ \mathbf{p} - \frac{1}{\kappa}_1 \mathbf{n} \\ \frac{1}{\kappa}_1 \\
	\end{pmatrix} \mapsto
	\begin{pmatrix}
		1 \\ \trans{\mathbf{p}} - \frac{1}{\trans{\kappa}_1} \frac{\n}{\|\n\|} \\ \frac{1}{\trans{\kappa}_1} \\
	\end{pmatrix} ,\quad
	\begin{pmatrix}
		1 \\ \mathbf{p} - \frac{1}{\kappa}_2 \mathbf{n} \\ \frac{1}{\kappa}_2 \\
	\end{pmatrix}\mapsto
	\begin{pmatrix}
		1 \\ \trans{\mathbf{p}} - \frac{1}{\trans{\kappa}_2} \frac{\n}{\|\n\|} \\ \frac{1}{\trans{\kappa}_2} \\
	\end{pmatrix} \quad
	\label{laguerre_tri_trans}
\end{equation}
For each transformation in Eq.~\eqref{laguerre_tri_trans}, the transformation in Eq.~\eqref{Laguerre}  is applied numerically to the left-hand-side vecor, and compared the result with the right-hand-side vector to compute $\|\n\|$, $\trans{\kappa}_1$, and $ \trans{\kappa}_2$. An L-isothermic surface satisfies the relation
\begin{equation}
	\trans{\kappa_1}\trans{A_1} = \trans{\kappa_2}\trans{A_2}
\end{equation}
Furthermore, the following relation is derived for evaluating the third fundamental form after transformation:
\begin{equation}
	\trans{\kappa}_1\trans{A}_1 = \frac{1}{\|\n\|}\kappa_1A_1 = \frac{1}{\lambda_n}e^\theta
\end{equation}
Finally we obtain the expressions of $\trans{A}_1$ and $\trans{A}_2$ as
\begin{equation}
	\trans{A_1} =  \frac{1}{\|\n\| \trans{\kappa}_1}e^\theta,\quad \trans{A_2} =  \frac{1}{\|\n\| \trans{\kappa}_2}e^\theta
\end{equation}
The stress distribution is obtained by incorporating $\trans{\kappa}_1$, $\trans{A}_1$, $\trans{\kappa}_2$, and $\trans{A}_2$ into Eq.~\eqref{trans_tention_equ}.
\end{document}